\documentclass[12pt]{article}
\usepackage[centertags]{amsmath}
\usepackage{amsfonts,amssymb,amsthm}

\newtheorem{thm}{Theorem}[section]

\newtheorem{lem}[thm]{Lemma}

\numberwithin{equation}{section}

\author{Tong Zhu}
\title{A Probability Method to Prove Combinatorial Identities}
\date{October 2009}

\begin{document}

\maketitle

\begin{abstract}
A probability method is provided to prove three classes of combinatorial identities. The method is extremely simple, only one step after the proper probability setup. 
\end{abstract}

\section{Introduction}
The probability method has been a popular and powerful tool to solve discrete problems in the last half century since Paul Erd\H os. It creates a new branch of mathematical research \cite{AS}. One of the beauties of this method is that there is being creative and open-minded involved, especially when the discrete structure and the probabilistic construction look totally irrelevant at first sight. This paper gives a probability method to prove some combinatorial identities.\\

The first result is a class of combinatorial identities involving symmetric functions,
\begin{align}
\sum_{j=1}^{n} (x_{j})^{m} \prod_{i \neq j} (x_{j}-x_{i})^{-1} = h_{m-n+1}(x),         \label{4}   
\end{align}
where $m \in \mathbb{Z}^{+} \cup \{0\}$ and $h_{m-n+1}(x)$ is the $(m-n+1)$-th homogeneous symmetric function of $x_{1}, ..., x_{n}$. \\

The combinatorial methods were provided by Louck and Biedenharn \cite{LB}, Strehl and Wilf \cite{SW} and Chen and Louck \cite{CL}. Equation \eqref{4} is a generalization of Good's identity \cite{Go},
\begin{align}
\sum_{j=1}^{n} \prod_{i \neq j} \Big(1- \frac{x_{j}}{x_{i}} \Big) ^{-1}=1.                            \label{3} 
\end{align}

The summation on the left hand side of \eqref{4} is important in the explicit matrix elements in the unitary groups, as well as multivariable hypergeometric series well poised in $SU(n)$. More combinatorial and algebraic applications of \eqref{4} and \eqref{3} can be found in \cite{GM} and \cite{LB}. \\

The second result contains multinomial convolutions.\\

\begin{align*}
&\prod_{i=1}^{n} \Big(\frac{\beta_{1}}{\beta_{i}} \Big)^{\alpha_{i}} \beta_{1}^{m}  \sum_{j=0}^{\infty} \delta_{j} \frac{(A+j+m-1)!}{(A+j-1)!} \\
=&
\sum_{i_{1}+ \cdots +i_{n}=m} \binom{m}{i_{1} \cdots i_{n}} \frac{(\alpha_{1}+i_{1}-1)!}{(\alpha_{1}-1)!} \beta_{1}^{i_{1}} \cdots \frac{(\alpha_{n}+i_{n}-1)!}{(\alpha_{n}-1)!} \beta_{n}^{i_{n}},
\end{align*}
where $m \geq 1$, $\alpha_{i}>0$, $\beta_{i}>0$ and $A=\sum_{i=1}^{n} \alpha_{i}$. The $\{\delta_{j}\}$ are defined recursively in Section 3.\\

A special case is about the convolution of two sequences.

\begin{align*}
\prod_{i=1}^{n} \Big(\frac{\beta_{1}}{\beta_{i}} \Big)^{\alpha_{i}} \beta_{1} \sum_{j=0}^{\infty} \delta_{j} (A+j)=\sum_{i=1}^{n} \alpha_{i}\beta_{i},
\end{align*}
The combinatorial proof of this type of identities was discussed in \cite{Ro} and \cite{Ha}.\\

The third result gives an identity involving multinomial convolution too.

\begin{align*}
 \frac{1}{ (n-1)! \displaystyle \prod_{i=1}^{n} a_{i}} \Big[ \frac{A^{m+n}}{m+n}+ B(n) \Big]
=\sum_{i_{1}+ \cdots + i_{n}=m} \binom{m}{i_{1}, \cdots, i_{n}} \frac{a_{1}^{i_{1}}}{i_{1}+1} \cdots \frac{a_{n}^{i_{n}}}{i_{n}+1}
\end{align*}
where $A=\sum_{i=1}^{n} a_{i}$ and
\begin{align*}
B(n)= \sum_{i=1}^{n} (-1)^{i}   \displaystyle \sum_{1 \leq j_{1} < j_{2} < \cdots <j_{i} \leq n}  \Big[  \sum_{k=0}^{n-1} \binom{n-1}{k} \Big( - \sum_{l=1}^{i} a_{j_{l}} \Big)^{n-1-k} \frac{ A^{k+m+1} - \big(\displaystyle \sum_{l=1}^{i} a_{j_{l}} \big)^{k+m+1}}{k+m+1}           \Big]       \Big\}   .
\end{align*}

A similar equation of the same type involving Bernoulli numbers was proved by Dilcher \cite{Di}.\\

A special case provides another method to prove an identity about the Stirling numbers of the second kind.

\begin{align*}
S(m+n,n)=\frac{1}{n!} \sum_{i=0}^{n} (-1)^{n-i} \binom{n}{i} i^{m+n}
\end{align*}

This probability method takes advantage of the nice structure of some probability density functions. In the next three sections, the combinatorial identities are proved by the summation of exponential random variables, Gamma random variables and uniform random variables. The last section discusses some inherent connections between the probability density functions and the combinatorial structures. 

\section{Exponential Distribution}

The probability density function of an exponential distribution with parameter $\lambda>0$ is
\begin{align*}
f(x)= \left \{ \begin{array}{lr}  \lambda e^{-\lambda x} & \text{if } x \geq 0 \\
0 & \text{if } x <0
\end{array} \right.
\end{align*}

In this paper, an important component of the proof is the density function of the summation of independent random variables of some specific distribution. The probability density function of the sum of independent exponential random variables with distinct parameters is as follows, which can be found in Feller ~\cite{Fe}, \#12, page 40 (with no proof). 

\vspace{0.5cm}

\begin{lem}
Let $\{X_{1},...,X_{n}\}$ $(n\geq 2)$ be a sequence of independent exponentially distributed random
variables with distinct parameters $\lambda_{1},...,\lambda_{n}$ $(\lambda_{i} >0)$, respectively. The density
function of their summation $\sum_{j=1}^{n} X_{j}$ is
\begin{align}
f_{n}(x)=\Big[\prod_{j=1}^{n}\lambda_{j} \Big]
\sum_{k=1}^{n}\frac{e^{-\lambda_{k}x}}{\displaystyle \prod_{\substack{l=1 \\
l \neq k}}^{n}(\lambda_{l}-\lambda_{k})}, \qquad x \geq 0.                                                             \label{2}
\end{align}
\end{lem}

\textit{Proof}: When $n=2$, by convolution
\begin{displaymath}
f_{2}(x)=\lambda_{1}\lambda_{2}
\frac{e^{-\lambda_{1}x}-e^{-\lambda_{2}x}}{\lambda_{2}-\lambda_{1}}.
\end{displaymath}
Suppose it is true for $n-1$. Then from convolution, we can get
\begin{displaymath}
f_{n}(x)=\Big[\prod_{j=1}^{n}\lambda_{j}
\Big]\sum_{k=1}^{n-1}\frac{e^{-\lambda_{k}x}-e^{-\lambda_{n}x}}{\displaystyle
\prod_{\substack{l=1 \\ l \neq k}}^{n}(\lambda_{l}-\lambda_{k})}
\end{displaymath}
In order to finish the proof, we need to show that
\begin{displaymath}
-\sum_{k=1}^{n-1} \frac{1}{\displaystyle \prod_{\substack{l=1 \\ l
\neq k}}^{n}(\lambda_{l}-\lambda_{k})}=\frac{1}{\displaystyle
\prod_{l=1}^{n-1}(\lambda_{l}-\lambda_{n})}
\end{displaymath}
Or equivalently,
\begin{displaymath}
\sum_{k=1}^{n} \frac{1}{\displaystyle \prod_{\substack{l=1 \\ l \neq
k}}^{n}(\lambda_{l}-\lambda_{k})}=0
\end{displaymath}
Multiplying the common denominator on both sides, this is equivalent to
\begin{align}
\sum_{k=1}^{n} (-1)^{k} \displaystyle \prod_{\substack{1\leq j<l\leq
n \\ l \neq k \\ j \neq k}}(\lambda_{l}-\lambda_{j})=0                                                                     \label{1}
\end{align}
Note that the product is a Vandermonde determinant:
\begin{displaymath}
\left| \begin{array}{ccccc} 1 & \lambda_{1} & \lambda_{1}^{2} & ...
&
\lambda_{1}^{n-2} \\
... & ... & ... & ... &...\\ 1 & \lambda_{k-1} & \lambda_{k-1}^{2} &
... &
\lambda_{k-1}^{n-2} \\
1 & \lambda_{k+1} & \lambda_{k+1}^{2} & ... &
\lambda_{k+1}^{n-2} \\
... & ...& ... & ... &...\\
1 & \lambda_{n} & \lambda_{n}^{2} & ... & \lambda_{n}^{n-2}
\end{array} \right|
\end{displaymath}
From this, it is not hard to see that \eqref{1} is the determinant of
\begin{displaymath}
\left( \begin{array}{ccccc} 1 & 1 & \lambda_{1} &  ... &
\lambda_{1}^{n-2} \\
... & ... & ... & ... &...\\ 1 & 1 & \lambda_{k} & ... &
\lambda_{k}^{n-2} \\
... & ...& ... & ... &...\\
1 & 1 & \lambda_{n} &  ... & \lambda_{n}^{n-2}
\end{array} \right),
\end{displaymath}
which is obviously zero. $\blacksquare$\\

\begin{thm}
For any distinct positive real values $\lambda_{1},...,\lambda_{n}$ ($n \geq 2$) and nonnegative integer $m$, 
\begin{equation*}
\sum_{k=1}^{n} \frac{1}{\lambda_{k}^{m}} \prod_{l \neq k} \frac{\lambda_{l}}{\lambda_{l}-\lambda_{k}}=\displaystyle \sum_{\substack{i_{1}+...+i_{n}=m \\ i_{1},...,i_{n} \geq 0}} \frac{1}{\lambda_{1}^{i_{1}}...\lambda_{n}^{i_{n}}}
\end{equation*}
\end{thm}

\textit{Proof} Let $\{X_{j}\}$, $1 \leq j  \leq n$, be independent exponential random variables with parameter $\lambda_{j}$, respectively. The $m$-th moment of $\sum X_{k}$ is
\begin{align*}
\mathbb{E}\Big[ (\sum_{k=1}^{n} X_{k})^{m} \Big]&=\mathbb{E} \Big[ \sum {m \choose i_{1}\; ...\; i_{n}} X_{1}^{i_{1}}...X_{n}^{i_{n}} \Big]\\
&=\sum_{i_{1}+...+i_{n}=m} \frac{m!}{\lambda_{1}^{i_{1}}...\lambda_{n}^{i_{n}}}
\end{align*}

On the other hand, by integrating the density function \eqref{2}, we have
\begin{equation*}
\mathbb{E}\Big[ (\sum_{k=1}^{n} X_{k})^{m} \Big]=\sum_{k=1}^{n} \frac{m!}{\lambda_{k}^{m}} \prod_{l \neq k} \frac{\lambda_{l}}{\lambda_{l}-\lambda_{k}}
\end{equation*}
$\blacksquare$\\

Equation \eqref{4} is an easy consequence of Theorem 2.2. \\



\section{Gamma Distribution}

A Gamma random variable $X$ has probability density function

\begin{align*}
f(x)=x^{\alpha-1} \frac{e^{-x/\beta}}{\Gamma(\alpha) \beta^{\alpha}}, \qquad x\geq 0,
\end{align*}
where  $\alpha>0$ and $\beta>0$ are two parameters.\\

The density function of the summation of independent Gamma random variables is given in \cite{Mo}.

\begin{lem} \label{gamma}
Let $X_{1}, \cdots,  X_{n}$ be independent Gamma random variables. The probability density function of $X_{i}$ is
\begin{align*}
f_{i}(x)=x^{\alpha_{i}-1} \frac{e^{-x/\beta_{i}}}{\Gamma(\alpha_{i}) \beta_{i}^{\alpha_{i}}}, \qquad x\geq 0
\end{align*}
with $\alpha_{i} >0$ and $\beta_{i}>0$. Without loss of generality, we may assume that $\beta_{1}=\displaystyle \min_{1 \leq i \leq n} \{\beta_{i}\}$. \\

Let $S_{n}=\sum_{i=1}^{n} X_{i}$. Then the density function of $S_{n}$ is
\begin{align*}
g_{n}(x)=\prod_{i=1}^{n} \Big(\frac{\beta_{1}}{\beta_{i}} \Big)^{\alpha_{i}} \sum_{j=0}^{\infty} \delta_{j} x^{A+j-1} \frac{e^{-x/\beta_{1}}}{\Gamma(A+j) \beta_{1}^{A+j}}, \qquad x \geq 0
\end{align*}
where $A=\sum_{i=1}^{n} \alpha_{i}$ and $\delta_{j}$'s are defined by the recursion
\begin{align*}
\delta_{j+1}=\frac{1}{j+1} \sum_{l=1}^{j+1} l \gamma_{l} \delta_{j+1-l}, \qquad j=0, 1, 2, \cdots
\end{align*}
with $\delta_{0}=1$ and
\begin{align*}
\gamma_{l}=\sum_{i=1}^{n} \alpha_{i} (1-\beta_{1}/\beta_{i})^{l}/l, \qquad l=1, 2, \cdots
\end{align*}
\end{lem}

Using two methods to compute the expectation of $S_{n}$, we can get an identity involving the convolution of the two sequences $\{\alpha_{i} \}$ and $\{\beta_{i}\}$.

\begin{thm}
\begin{align*}
\prod_{i=1}^{n} \Big(\frac{\beta_{1}}{\beta_{i}} \Big)^{\alpha_{i}} \beta_{1} \sum_{j=0}^{\infty} \delta_{j} (A+j)=\sum_{i=1}^{n} \alpha_{i}\beta_{i},
\end{align*}
with $A$ and $\delta_{j}$ as defined in Lemma \ref{gamma}.
\end{thm}

From the computation of moments, a general identity is derived.

\begin{thm}
If $m$ is a positive integer,
\begin{align*}
&\prod_{i=1}^{n} \Big(\frac{\beta_{1}}{\beta_{i}} \Big)^{\alpha_{i}} \beta_{1}^{m}  \sum_{j=0}^{\infty} \delta_{j} \frac{(A+j+m-1)!}{(A+j-1)!} \\
=&
\sum_{i_{1}+ \cdots +i_{n}=m} \binom{m}{i_{1} \cdots i_{n}} \frac{(\alpha_{1}+i_{1}-1)!}{(\alpha_{1}-1)!} \beta_{1}^{i_{1}} \cdots \frac{(\alpha_{n}+i_{n}-1)!}{(\alpha_{n}-1)!} \beta_{n}^{i_{n}}
\end{align*}
\end{thm}


\section{Uniform Distribution}

A uniformly distributed random variable $X$ on the interval $[a, b]$ has probability density function
\begin{displaymath}
f(x)= \left \{ 
\begin{array}{rl}
\frac{1}{b-a} & \text{if } a \leq x \leq b \\

0  & \text{otherwise}
\end{array} \right.
\end{displaymath}

Two types of summations are considered. One is the summation of I.I.D. (independent and identically distributed) uniform random variables. The other one is the summation of independent random variables uniformly distributed on different intervals. Before the discussion, the following notation is needed.
\begin{displaymath}
x_{+}=\frac{x+|x|}{2}, \text{ for } x \in \mathbb{R}.
\end{displaymath}

Let $X_{1}, ..., X_{n}$ be uniformly distributed I.I.D. random variables. Without loss of generality, we may assume that $X_{i}$ is uniform on $[0,a]$. (If $X_{i}$ is uniformly distributed on $[b,c]$, then let $Y_{i}=X_{i}-b$. So $Y_{i}$'s are I.I.D. random variables uniformly distributed on $[0, c-b]$ and $\sum_{i=1}^{n} X_{i}=\sum_{i=1}^{n} Y_{i}+nb$. And it is sufficient to consider $\sum_{i=1}^{n} Y_{i}$.) The following lemma, which can be found in \cite{Fe} page 27,  gives the density function of the summation.

\begin{lem}
For $n\geq 1$, let $S_{n}=\sum_{i=1}^{n} X_{i}$. The probability density function of $S_{n}$ is
\begin{align}
h_{n}(x)=\frac{1}{a^{n}(n-1)!} \sum_{i=0}^{n} (-1)^{i} \binom{n}{i} (x-ia)_{+}^{n-1},    \label{uniform density 1}
\end{align}
for $x \in [0, na]$ and $h_{n}(x)=0$ otherwise.
\end{lem}

Using two methods to compute the $m$-th moment of $S_{n}$, we can get

\begin{lem}
For $m \geq 1$, 
\begin{align}
&\frac{1}{(n-1)!}  \Big\{ \frac{n^{m+n}}{m+n}+    \sum_{i=1}^{n} (-1)^{i} \binom{n}{i}  \Big[ \sum_{j=0}^{n-1} \binom{n-1}{j} (-i)^{n-j-1} \frac{n^{m+j+1}-i^{m+j+1}}{m+j+1} \Big]   \Big\} \nonumber \\
=& \sum_{i_{1}+ \cdots + i_{n}=m} \binom{m}{i_{1}, \cdots, i_{n}} \frac{1}{(i_{1}+1) \cdots (i_{n}+1)}   \label{uniform identity 1}
\end{align}
\end{lem}

The right side of the above identity can be simplified. 

\begin{lem}
\begin{align*}
\sum_{i_{1}+ \cdots + i_{n}=m} \binom{m}{i_{1}, \cdots, i_{n}} \frac{1}{(i_{1}+1) \cdots (i_{n}+1)}=\frac{S(m+n,n)}{\binom{m+n}{n}},   
\end{align*}
where $S(m+n,n)$ is the Stirling number of the second kind. 
\end{lem}

\textit{Proof}: Let
\begin{align*}
\gamma_{m,n}=\sum_{i_{1}+ \cdots + i_{n}=m} \binom{m}{i_{1}, \cdots, i_{n}} \frac{1}{(i_{1}+1) \cdots (i_{n}+1)}.
\end{align*}

Note that $\gamma_{m,n}$ is the coefficient of $t^{m+n}/m!$ in $(e^{t}-1)^{n}$. Hence,
\begin{align*}
e^{y(e^{t}-1)}&=\sum_{n \geq 0} \frac{y^{n}}{n!} (e^{t}-1)^{n} \\
&=\sum_{n \geq 0} \frac{y^{n}}{n!} \sum_{m \geq 0} \frac{t^{m+n}}{m!} \gamma_{m,n}
\end{align*}

Besides,
\begin{align*}
e^{y(e^{t}-1)}&=\sum_{n\geq 0} y^{n} \Big[ \sum_{k \geq n} S(k,n) \frac{t^{k}}{k!} \Big] \\
&=\sum_{n \geq 0} y^{n} \Big[\sum_{m \geq 0} S(m+n,n) \frac{t^{m+n}}{(m+n)!} \Big]
\end{align*}

Comparing the coefficients, there is
\begin{align*}
\frac{1}{n!} \sum_{m \geq 0} \frac{t^{m+n}}{m!} \gamma_{m,n} = \sum_{m \geq 0} S(m+n,n) \frac{t^{m+n}}{(m+n)!}.
\end{align*}

So
\begin{align*}
\frac{\gamma_{m,n}}{n!}=S(m+n,n) \frac{m!}{(m+n)!}.
\end{align*}
$\blacksquare$\\

Observe the density function \eqref{uniform density 1}. The $m$-th moment of $S_{n}$ can be represented in another way, instead of the two sides in \eqref{uniform identity 1}.

\begin{lem}
\begin{align*}
\frac{1}{(n-1)!} \sum_{i=0}^{n} (-1)^{i} \binom{n}{i} \int_{0}^{n} (x-i)_{+}^{n-1} x^{m} dx =\frac{m!}{(m+n)!} \sum_{i=0}^{n} (-1)^{n-i} \binom{n}{i} i^{m+n}
\end{align*}
\end{lem}

\textit{Proof}: First,
\begin{align}
 \int_{0}^{n} (x-i)_{+}^{n-1} x^{m} dx &= \int_{i}^{n} (x-i)^{n-1} x^{m} dx \nonumber\\
 &=\frac{m!(n-1)!}{(m+n)!} \sum_{j=0}^{m} (-1)^{j} \binom{m+n}{m-j} (n-i)^{n+j} n^{m-j} \nonumber\\
&= \frac{m!(n-1)!}{(m+n)!}(-1)^{n} \sum_{j=0}^{m} (-1)^{n+j} \binom{m+n}{m-j} (n-i)^{n+j} n^{m-j} \label{uni step 1}
\end{align}

Note that the summation in \eqref{uni step 1} is part of the expansion of $i^{m+n}$.
\begin{align*}
i^{m+n}=\sum_{j=0}^{m+n} (-1)^{j} \binom{m+n}{j} (n-i)^{j} n^{m+n-j}
\end{align*}

But by symmetry, the other part will vanish after the summation of $i$.
\begin{align*}
\sum_{i=0}^{n} (-1)^{n-i} \binom{n}{i} \Big[ \sum_{j=0}^{n-1} (-1)^{j} \binom{m+n}{j} (n-i)^{j} n^{m+n-j}  \Big]=0
\end{align*}
$\blacksquare$ \\

Hence, the uniform density function gives another method to prove the identity of Stirling function of the second kind.

\begin{thm}
\begin{align*}
S(m+n,n)=\frac{1}{n!} \sum_{i=0}^{n} (-1)^{n-i} \binom{n}{i} i^{m+n}
\end{align*}
\end{thm}

\vspace{1cm}

If  $X_{1}, ..., X_{n}$ are independent random variables and $X_{i}$ is uniformly distributed on $[0, a_{i}]$, the density function of their summation is given by \cite{SNH}.

\begin{lem}
For $n \geq 2$, let $S_{n}=\sum_{i=1}^{n} X_{i}$. The probability density function of $S_{n}$ is
\begin{align*}
\tilde{h}_{n}(x)= \frac{1}{ (n-1)! \displaystyle \prod_{i=1}^{n} a_{i}} \Big\{ x^{n-1} + \sum_{i=1}^{n} (-1)^{i} \displaystyle \sum_{1 \leq j_{1} < j_{2} < \cdots <j_{i} \leq n} 
\Big[ \Big( x-\sum_{l=1}^{i} a_{j_{l}}       \Big)_{+}  \Big] ^{n-1}      \Big\}.
\end{align*}
\end{lem}

\begin{thm}
For $m \geq 1$, 
\begin{align*}
 \frac{1}{ (n-1)! \displaystyle \prod_{i=1}^{n} a_{i}} \Big[ \frac{A^{m+n}}{m+n}+ B(n) \Big]
=\sum_{i_{1}+ \cdots + i_{n}=m} \binom{m}{i_{1}, \cdots, i_{n}} \frac{a_{1}^{i_{1}}}{i_{1}+1} \cdots \frac{a_{n}^{i_{n}}}{i_{n}+1}
\end{align*}
where $A=\sum_{i=1}^{n} a_{i}$ and
\begin{align*}
B(n)= \sum_{i=1}^{n} (-1)^{i}   \displaystyle \sum_{1 \leq j_{1} < j_{2} < \cdots <j_{i} \leq n}  \Big[  \sum_{k=0}^{n-1} \binom{n-1}{k} \Big( - \sum_{l=1}^{i} a_{j_{l}} \Big)^{n-1-k} \frac{ A^{k+m+1} - \big(\displaystyle \sum_{l=1}^{i} a_{j_{l}} \big)^{k+m+1}}{k+m+1}           \Big]       \Big\}   .
\end{align*}
\end{thm}



\section{Relations and Future Development}

It is easy to start with the summation of some independent random variables to get some identities, if the density function is provided. But given an identity, how to construct a probability distribution for the proof is another question. Hence, it is important to find out the inherent relations between the probability density functions and the combinatorial structures. One possible connection is their ``generating functions''. All of the identities in this paper were studied by the generating function method before. The generating function of a probability distribution is its \textit{characteristic function} (ch.f.). The ch.f. of a random variable is defined as
\begin{align*}
\phi_{X}(t)=\mathbb{E}(e^{itX}).
\end{align*}

For example, the ch.f. of an exponential distribution with parameter $\lambda$ is
\begin{align*}
\phi(t)=\frac{1}{1-\frac{it}{\lambda}}.
\end{align*}

And Lemma 2.1 can also be proved by the ch.f. $\phi_{n}(t)$ of the summation $\displaystyle \sum_{j=1}^{n} X_{j}$
\begin{align}
\phi_{n}(t)&=\prod_{j=1}^{n} \mathbb{E} e^{itX_{j}} \nonumber\\
&=\prod_{j=1}^{n} \frac{1}{1-\frac{it}{\lambda_{j}}}.                                \label{5}
\end{align}

Note that this characteristic function has the same form as the generating function of the homogeneous symmetric functions $h_{r}(x)$ in variables $x_{1}, x_{2}, ...$
\begin{align*}
G(t)=\sum_{r \geq 0} h_{r} t^{r}=\prod_{j \geq 1} (1-x_{j}t)^{-1}.
\end{align*}

Besides, the density function \eqref{2} has similar product  inside the summation to the form of Lagrange interpolation polynomial
\begin{align*}
p(x)=\sum_{j=1}^{n} p(x_{j}) \prod_{k \neq j} \frac{x-x_{k}}{x_{j}-x_{k}}.
\end{align*}

In fact, Good's identity \eqref{3} and its generalization \eqref{4} can also be proved by Lagrange interpolation method (as in \cite{CL} and \cite{SW}).\\

If we compute the expectation $\mathbb{E}(\exp(t\sum X_{k}))$ by integrating the integral transform of \eqref{2} and by the fact that $\mathbb{E}(\exp(t\sum X_{k}))=\prod \mathbb{E}(\exp(t X_{k}))$, we get another identity.
\begin{align}
\sum_{k=1}^{n} \frac{1}{1-\frac{t}{\lambda_{k}}} \prod_{l \neq k} \frac{\lambda_{l}}{\lambda_{l}-\lambda_{k}} = \prod_{k=1}^{n} \frac{1}{1-\frac{t}{\lambda_{k}}},\qquad t<\min\{\lambda_{1}, ..., \lambda_{n}\}.                                                              \label{6}                  
\end{align}

Observe that \eqref{6} can deduce the interpolation proof of Theorem 1 in \cite{CL} and \cite{SW}. \\

\vspace{0.4cm} 

We may also check the  \textit{moment generating function} (m.g.f.) of a distribution, if it exists. The m.g.f. of a random variable $X$ is defined as 
\begin{align*}
\mathcal{M}_{X}(t)=\mathbb{E}(e^{tX}).
\end{align*}

The m.g.f. also has the property that the m.g.f. of the summation of independent random variables is the product of their m.g.f.'s. And the moments can be got from the m.g.f..
\begin{align*}
\mathbb{E}(X^{m})=\mathcal{M}_{X}^{(m)}(0).
\end{align*}

For example, the m.g.f. of the uniform distribution on the interval $[0,1]$ is
\begin{align*}
\mathcal{M}(t)=\frac{e^{t}-1}{t}.
\end{align*}

Then the m.g.f. of $S_{n}$, the summation of $n$ independent uniform random variables on $[0,1]$, is
\begin{align*}
\mathcal{M}_{n}(t)=\frac{(e^{t}-1)^{n}}{t^{n}}.
\end{align*}

Note that the exponential generating function of the Stirling number of the second kind has similar form.
\begin{align*}
F_{n}(t)&=\sum_{m \geq 0} S(m+n)\frac{t^{m+n}}{(m+n)!}\\
&=\frac{(e^{t}-1)^{n}}{n!}.
\end{align*}

Look at the exponential generating function of $\binom{m+n}{n} \mathbb{E}(S_{n}^{m})$.
\begin{align*}
&\sum_{m \geq 0} \binom{m+n}{n} \mathbb{E}(S_{n}^{m}) \frac{t^{m+n}}{(m+n)!}= \frac{t^{n}}{n!} \sum_{m\geq 0} \frac{t^{m} \mathbb{E}(S_{n}^{m}) }{m!} \\
&=\frac{t^{n}}{n!} \frac{(e^{t}-1)^{n}}{t^{n}}=\frac{(e^{t}-1)^{n}}{n!}.
\end{align*}

This proofs Lemma 4.3.\\



\vspace{0.5cm}

\textbf{Acknowledgements.} I would like to thank Robin Pemantle and Herbert Wilf for the thoughtful comments and suggestions. Many thanks to James Haglund, Richard Stanley, Curtis Greene, Richard Askey and Stephen Milne for their kind references.

\end{document}